\documentclass[11pt,letterpaper,reqno]{amsart}

\title{Contact Geometry of the Pontryagin Maximum Principle}
\author{Tomoki Ohsawa}
\address{Department of Mathematical Sciences, The University of Texas at Dallas, 800 W Campbell Rd, Richardson, TX 75080-3021}
\email{tomoki@utdallas.edu}
\date{\today}

\keywords{Optimal control; Maximum principle; Differential geometric methods; Geometric approaches.}

\subjclass[2010]{37J55, 49J15, 53D10}

\usepackage{graphicx}
\usepackage[margin=1in, marginpar=.5in]{geometry}
\usepackage[numbers,sort&compress]{natbib}
\usepackage[all]{xy}
\usepackage{tikz-cd}

\usepackage[colorlinks=false]{hyperref}

\theoremstyle{plain}
\newtheorem{theorem}{Theorem}[section]

\newtheorem{proposition}[theorem]{Proposition}

\theoremstyle{definition}

\theoremstyle{remark}
\newtheorem{remark}[theorem]{Remark}

\def\od#1#2{\dfrac{d#1}{d#2}}
\def\pd#1#2{\dfrac{\partial #1}{\partial #2}}

\def\parentheses#1{\!\left(#1\right)}
\def\brackets#1{\!\left[#1\right]}
\def\braces#1{\!\left\{#1\right\}}

\def\DS{\displaystyle}
\def\R{\mathbb{R}}
\def\defeq{\mathrel{\mathop:}=}
\def\eqdef{=\mathrel{\mathop:}}
\def\setdef#1#2{ \left\{ #1 \ |\ #2 \right\} }
\def\ip#1#2{\left\langle#1,#2\right\rangle}
\def\tip#1#2{\langle#1,#2\rangle}

\def\d{{\bf d}}
\def\ins#1{{\bf i}_{#1}}

\begin{document}

\footskip=.6in

\begin{abstract}
  This paper gives a brief contact-geometric account of the Pontryagin maximum principle.
  We show that key notions in the Pontryagin maximum principle---such as the separating hyperplanes, costate, necessary condition, and normal/abnormal minimizers---have natural contact-geometric interpretations.
  We then exploit the contact-geometric formulation to give a simple derivation of the transversality condition for optimal control with terminal cost.
\end{abstract}

\maketitle

\section{Introduction}
It is well known that a necessary condition for optimality of the Pontryagin maximum principle may be interpreted as a Hamiltonian system, and so its geometric formulation usually exploits the language of symplectic geometry; see e.g., \citet[Chapter~11]{Ju1997}, \citet[Chapter~12]{AgSa2004}.

The main focus of this paper is to change this perspective slightly to look at the maximum principle from the point of view of contact geometry, the ``odd-dimensional cousin''~(\citet[Appendix~4]{Ar1991}) of symplectic geometry.
The correspondence between contact and symplectic Hamiltonian systems is elementary and well known (see, e.g., \citet[Appendix~4]{Ar1991}), and thus switching between symplectic and contact views is fairly trivial as far as the mathematical technicality is concerned.
Our stress here is rather that {\em the language of contact geometry fits more naturally to a proof of the Pontryagin maximum principle and so we may exploit the contact-geometric view from the outset.
It also provides an alternative geometric perspective on applications of the maximum principle}.

\section{Geometry of Optimal Control on $\R^{n}$}
\subsection{Extended System}
Let $\R^{n}$ be the state space and $\mathcal{U}$ be a compact subset of $\R^{m}$ that defines the space of controls.
Define a control system by $f\colon \R^{n} \times \mathcal{U} \to \R^{n}$, and let $L\colon \R^{n} \times \mathcal{U} \to \R$ be the cost function and $S_{1}$ be a submanifold of $\R^{n}$.
Then consider the following optimal control problem:
\begin{equation}
  \label{eq:OptCtrlProblem}
  \begin{array}{c}
    \DS \min_{u(\cdot) \in \mathcal{U}} \int_{t_{0}}^{t_{1}} L(x(t), u(t))\,dt
    \medskip\\
    \DS \text{subject to}
    \quad
    \dot{x} = f(x, u),
    \quad
    x(t_{0}) = x_{0},
    \quad
    x(t_{1}) \in S_{1},
  \end{array}
\end{equation}
where the initial time $t_{0} \in \R$ and initial point $x_{0} \in \R^{n}$ are fixed whereas the terminal time $t_{1} \in \R$ is free.

Recall (see, e.g., \citet{PoBoGaMi1962}, \citet[Chapter~4]{Li2012} and \citet{Lewis-PMP}) that the first step in proving the Pontryagin maximum principle is to introduce a new variable (running cost) $x^{0}$ by
\begin{equation*}
  x^{0}(t) \defeq \int_{t_{0}}^{t} L(x(s), u(s))\,ds,
\end{equation*}
that is, $x^{0}$ may be regarded as a solution of the differential equation
\begin{equation*}
  \dot{x}^{0} = L(x, u)
\end{equation*}
with the initial condition $x^{0}(t_{0}) = 0$.
One then augments the original control system by the above system: We define the extended state variable
\begin{equation*}
  \hat{x} \defeq (x^{0}, x) \in \R^{n+1},
\end{equation*}
and define $\hat{f}\colon \R^{n+1} \times \mathcal{U} \to \R^{n+1}$ by
\begin{equation*}
  \hat{f}(\hat{x},u) \defeq 
  \begin{bmatrix}
    L(x,u) \\
    f(x,u)
  \end{bmatrix}.
\end{equation*}
Then the optimal control problem \eqref{eq:OptCtrlProblem} is restated as
\begin{equation*}
 \min_{u(\cdot) \in \mathcal{U}} x^{0}(t_{1})
\end{equation*}
subject to the {\em extended system} (sometimes called the Mayer form; see, e.g., \citet[Chapter~4]{Li2012}) on $\R^{n+1}$ defined by
\begin{equation}
  \label{eq:extended_system}
  \dot{\hat{x}} = \hat{f}(\hat{x}, u)
\end{equation}
along with the end points $\hat{x}(t_{0}) = (0, x_{0}) \eqdef \hat{x}_{0}$ and $\hat{x}(t_{1}) = (x^{0}(t_{1}), x_{1}) \eqdef \hat{x}_{1}$.

\subsection{Costate Lives in a Projective Space}
\label{ssec:Costate_and_ProjectiveSpace}
Now, let $u^{\star}\colon [t_{0}, t_{1}^{\star}] \to \mathcal{U}$ be an optimal control and $\hat{x}^{\star}\colon [t_{0}, t_{1}^{\star}] \to \R^{n+1}$ be the corresponding optimal trajectory of the extended system~\eqref{eq:extended_system}.
Combining needle variations and temporal variations at the terminal time $t_{1}^{\star}$ of the optimal control $u^{\star}$ gives the terminal cone $C_{\hat{x}_1^{\star}} \subset \R^{n+1}$ that approximates the reachable set near the terminal point $\hat{x}^{\star}_{1} = (x^{\star0}_{1}, x^{\star}_{1}) \defeq \hat{x}^{\star}(t_1^{\star})$ (see, e.g., \citet[Section~4.2]{Li2012} and \citet[Chapter~5]{Lewis-PMP} for details of the construction of the terminal cone).

\begin{figure}[htbp]
  \begin{minipage}{0.45\linewidth}
    \centering
    \includegraphics[width=.85\linewidth]{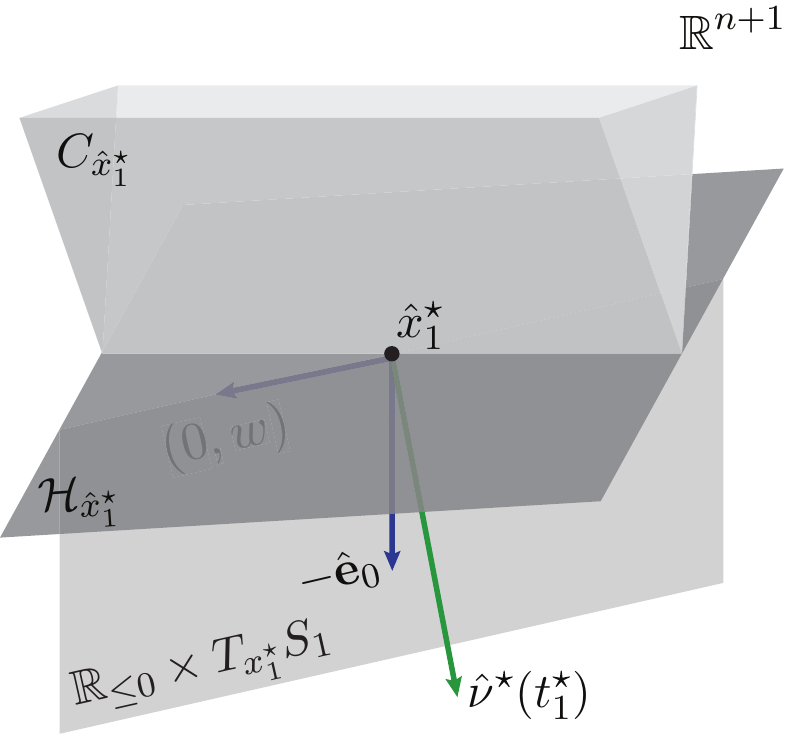}
    \caption{Terminal cone $C_{\hat{x}_1^{\star}}$, separating hyperplane $\mathcal{H}_{\hat{x}_1^{\star}}$, and costate $\hat{\nu}^{\star}(t_{1}^{\star})$.}
    \label{fig:SeparatingHyperplane}
  \end{minipage}
  \quad
  \begin{minipage}{0.45\linewidth}
  \centering
  \includegraphics[width=.55\linewidth]{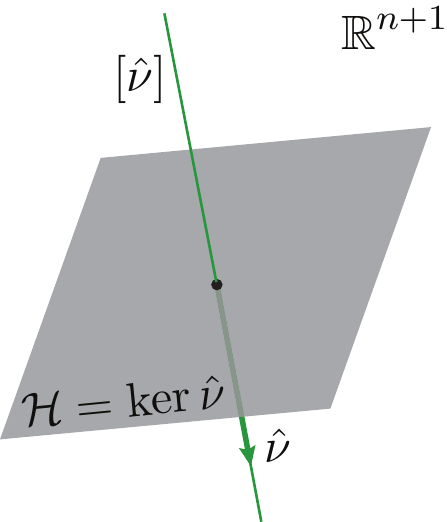}
  \caption{Hyperplane $\mathcal{H}$ and costate $[\hat{\nu}]$}
  \label{fig:HyperplaneAndRay}
  \end{minipage}
\end{figure}

One then argues that the interior of the cone $C_{\hat{x}_1^{\star}}$ does not intersect $\R_{\le0} \times T_{x^{\star}_1}S_{1}$, where $\R_{\le0}$ is the set of non-positive real numbers and $T_{x^{\star}_1}S_{1}$ is the tangent space to $S_{1}$ at $\hat{x}^{\star}_{1}$; because if it did then that implies that there exists a variation of the optimal control $u^{\star}$ with the terminal point still in $S_{1}$ but with a lower total cost; so $\R_{\le0} \times T_{x^{\star}_1}S_{1}$ defines ``forbidden'' directions.
As a result, one concludes that there exists a hyperplane $\mathcal{H}_{\hat{x}_1^{\star}} \subset \R^{n+1}$ that separates the interior of $C_{\hat{x}_1^{\star}}$ and $\R_{\le0} \times T_{x^{\star}_1}S_{1}$ in the sense that they sit on different sides from each other (see Fig.~\ref{fig:SeparatingHyperplane}).

One then introduces the costate vector $\hat{\nu}^{\star}(t_1^{\star}) \in (\R^{n+1})^{*}\backslash\{0\} \cong \R^{n+1}\backslash\{0\}$ as an element such that $\ker\hat{\nu}^{\star}(t_1^{\star}) = \mathcal{H}_{\hat{x}_1^{\star}}$.
However, we observe that $\hat{\nu}^{\star}(t_1^{\star})$ is not uniquely defined: We may multiply $\hat{\nu}^{\star}(t_1^{\star})$ by any $k \in \R\backslash\{0\}$ and have $\ker(k\,\hat{\nu}^{\star}(t_1^{\star})) = \mathcal{H}_{\hat{x}_1^{\star}}$, i.e., two costate vectors $\hat{\nu}_{1}$ and $\hat{\nu}_{2}$ in $\R^{n+1}$ are equivalent if one is a nonzero constant multiple of the other:
\begin{equation*}
  \hat{\nu}_{1} \sim \hat{\nu}_{2} \iff \text{$\hat{\nu}_{2} = k\,\hat{\nu}_{1}$ for some $k \in \R\backslash\{0\}$}.
\end{equation*}
This defines an equivalence relation and thus we may define the equivalence class $[\hat{\nu}]$ of an element $\hat{\nu} \in \R^{n+1}$ by
\begin{equation*}
  [\hat{\nu}] \defeq \setdef{ \hat{\mu} \in \R^{n+1}\backslash\{0\} }{ \text{$\hat{\mu} = k\,\hat{\nu}$ for some $k \in \R\backslash\{0\}$} }.
\end{equation*}
Geometrically, the equivalence class $[\hat{\nu}]$ corresponds to the straight line along the vector $\hat{\nu}$, and the collection of these equivalence classes (straight lines passing through the origin) $[\hat{\nu}]$ defines the projective space $\mathbb{P}(\R^{n+1}) \defeq (\R^{n+1}\backslash\{0\})/{\sim}$.
Therefore, the costate is better defined as the straight line along the vector $\hat{\nu}^{\star}(t_1^{\star})$ than the vector itself (see Fig.~\ref{fig:HyperplaneAndRay}), or more mathematically speaking, {\em the costate is most naturally defined as an element in the projective space $\mathbb{P}(\R^{n+1})$}.

If we need to choose a representative element $\hat{\nu} = (\nu_{0}, \nu)$ of $[\hat{\nu}]$, we choose {\em by convention} an element $\hat{\nu} \in \R^{n+1}$ such that $\ip{\hat{\nu}}{\hat{w}} \ge 0$ for any $\hat{w} \in \R_{\le0} \times T_{x^{\star}_1}S_{1}$, where $\ip{\cdot}{\cdot}$ stands for the standard pairing of vectors in $\R^{n+1}$; that is, $\hat{\nu}$ and $\R_{\le0} \times T_{x^{\star}_1}S_{1}$ are on the same side as shown in Fig.~\ref{fig:SeparatingHyperplane}.
In particular, choosing $\hat{w} = -\hat{\bf e}_{0} = (-1,0,\dots,0)$ gives $\nu_{0} \le 0$ and $\hat{w} = (0,w)$ with arbitrary $w \in T_{x^{\star}_1}S_{1}$ gives the transversality condition:
\begin{equation}
  \label{eq:transversality_condition}
  \nu^{\star}(t_{1}^{\star}) \in (T_{x^{\star}_{1}}S_{1})^{\perp}.
\end{equation}

\subsection{Adjoint Equation and Control Hamiltonian}
The costate vector $\hat{\nu}^{\star}(t_{1}^{\star})$ defined above encodes a necessary condition for optimality at the terminal time $t_{1}^{\star}$.
We then propagate the costate $\hat{\nu}^{\star}(t_{1}^{\star})$ back along the optimal trajectory $\hat{x}^{\star}(t)$ to formulate a necessary condition {\em along} the trajectory.

This is done by defining the costate vectors $\hat{\nu}^{\star}\colon [t_{0}, t_{1}^{\star}] \to \R^{n+1}$ as the solution to the adjoint equation 
\begin{equation}
  \label{eq:adjoint}
  \dot{\nu}^{\star}_{0} = 0,
  \qquad
  \dot{\nu}^{\star}_{i} = -\nu^{\star}_{j} \pd{f^{j}}{x^{i}}(x^{\star}, u^{\star}) - \nu^{\star}_{0}\pd{L}{x^{i}}(x^{\star}, u^{\star})
\end{equation}
with the terminal condition $\hat{\nu}^{\star}(t_{1}^{\star})$; where the indices $i$ and $j$ run from $1$ to $n$.
The motivation for doing so is that one can relate the costate $\hat{\nu}^{\star}(t)$ and perturbation $\delta\hat{x}(t)$ at $\hat{x}^{\star}(t)$ with those at the terminal point $\hat{x}^{\star}_{1}$ thanks to the conservation of the pairing of them, i.e.,
\begin{equation}
  \label{eq:pairing_conservation}
  \ip{ \hat{\nu}^{\star}(t) }{ \delta\hat{x}(t) } = \ip{ \hat{\nu}^{\star}(t_{1}^{\star}) }{ \delta\hat{x}(t_{1}^{\star}) }.
\end{equation}
\begin{remark}
  \label{rem:cotangent_lift}
  The above conservation is due to the fact that the propagation of perturbation $\delta\hat{x}(t) \mapsto \delta\hat{x}(t_{1}^{\star})$ is the tangent lift $T\phi_{t_{1}^{\star}-t}$ of the flow $\phi_{t_{1}^{\star}-t}\colon \hat{x}^{\star}(t) \mapsto \hat{x}^{\star}(t_{1}^{\star})$ defined by the optimal solution $\dot{\hat{x}}^{\star} = \hat{f}(\hat{x}^{\star}, u^{\star})$ whereas the time-reversed propagation of costate $\hat{\nu}^{\star}(t_{1}^{\star}) \mapsto \hat{\nu}^{\star}(t)$ is the cotangent lift $T^{*}\phi_{t_{1}^{\star}-t}$ of $\phi_{t_{1}^{\star}-t}$.
  In fact, the adjoint equation~\eqref{eq:adjoint} is nothing but the time derivative of the cotangent lift $T^{*}\phi_{-t}$ (see, e.g., \citet[Chapter~12]{AgSa2004}).
\end{remark}
Recall (see, e.g., \cite{PoBoGaMi1962,Li2012,Lewis-PMP}) also that setting $\delta\hat{x}(t_{1}^{\star})$ equal to the variation resulting from a needle variation ending at $t$ in \eqref{eq:pairing_conservation} yields the following essential result of the maximum principle:
For any $u \in \mathcal{U}$,
\begin{equation*}
  H_{\rm c}(\hat{x}^{\star}(t), \hat{\nu}^{\star}(t), u) \le H_{\rm c}(\hat{x}^{\star}(t), \hat{\nu}^{\star}(t), u^{\star}(t)),
\end{equation*}
with the control Hamiltonian $H_{\rm c}\colon \R^{n+1} \times \R^{n+1} \times \mathcal{U} \to \R$ defined by
\begin{equation}
  \label{eq:hatH}
  H_{\rm c}(\hat{x}, \hat{\nu}, u) \defeq \tip{\hat{\nu}}{\hat{f}(x,u)} = \nu \cdot f(x, u) + \nu_{0} L(x, u).
\end{equation}
Therefore, we may define the optimal Hamiltonian $H\colon \R^{n+1} \times \R^{n+1} \to \R$ as follows:
\begin{equation*}
  H(\hat{x}, \hat{\nu}) \defeq \max_{u \in \mathcal{U}} H_{\rm c}(\hat{x}, \hat{\nu}, u) = H_{\rm c}(\hat{x}, \hat{\nu}, u^{\star}(\hat{x},\hat{\nu})).
\end{equation*}

\subsection{Contact Control Hamiltonian and Normal \& Abnormal Extremals}
Notice that the control Hamiltonian~\eqref{eq:hatH} is homogeneous of degree 1 in the costate vector $\hat{\nu} = (\nu_{0}, \nu) \in \R^{n+1}$, i.e., for any $k \in \R\backslash\{0\}$,
\begin{equation*}
  H_{\rm c}(\hat{x}, k\,\hat{\nu}, u) = k\,H_{\rm c}(\hat{x}, \hat{\nu}, u),
\end{equation*}
and hence, taking the quotient by $\R\backslash\{0\}$, it projects to the {\em contact control Hamiltonian}
\begin{equation}
  \label{eq:h_c}
  h_{\rm c}\colon \R^{n+1} \times \mathbb{P}(\R^{n+1}) \times \mathcal{U} \to \R.
\end{equation}
Likewise, we define the {\em optimal contact Hamiltonian} $h\colon \R^{n+1} \times \mathbb{P}(\R^{n+1}) \to \R$ as follows:
\begin{equation}
  \label{eq:h}
  h(\hat{x}, [\hat{\nu}]) \defeq \max_{u \in \mathcal{U}} h_{\rm c}(\hat{x}, [\hat{\nu}], u).
\end{equation}
The above definitions of control Hamiltonians are more natural for us because, as shown in Section~\ref{ssec:Costate_and_ProjectiveSpace}, the costate essentially lives in the projective space $\mathbb{P}(\R^{n+1})$.

The contact control Hamiltonian $h_{\rm c}$ is well-defined globally for any costate in $\mathbb{P}(\R^{n+1})$ regardless of whether the extremal in question is normal or abnormal.
However, it turns out that abnormal extremals fall into the coordinate singularity of the natural coordinates for normal extremals:
For normal extremals, i.e., if $\nu_{0} \neq 0$, then we set $\lambda \defeq -\nu/\nu_{0} \in \R^{n}$.
Note that $\lambda$ is nothing but {\em homogeneous coordinates} for $\mathbb{P}(\R^{n+1})$: A common way of giving coordinates to $\mathbb{P}(\R^{n+1})$ is to identify $[\hat{\nu}] = [\nu_{0}: \nu_{1}: \dots : \nu_{n}] \in \mathbb{P}(\R^{n+1})$ with $[-\hat{\nu}/\nu_{0}]$ and write
\begin{equation}
  \label{eq:homogeneous_coordinates}
  [\hat{\nu}] = [-\hat{\nu}/\nu_{0}]
  = [(-1,\lambda)]
  = \brackets{ -1: \lambda_{1}: \dots : \lambda_{n} },
\end{equation}
where $\{ \lambda_{i} \defeq -\nu_{i}/\nu_{0} \}_{i=1}^{n}$ are coordinates for $[\hat{\nu}]$.
Recall from Section~\ref{ssec:Costate_and_ProjectiveSpace} that we have $\nu_{0} \le 0$ by convention and so $\nu_{0} < 0$ here: The negative sign in the definition of $\lambda$ makes the above representative element $(-1, \lambda) \in \R^{n+1}$ meet this convention.
As a result, we may define the contact control Hamiltonian~\eqref{eq:h_c} as follows:
\begin{equation*}
  h_{\rm c}(\hat{x}, [\hat{\nu}], u) = H_{\rm c}(\hat{x}, -\hat{\nu}/\nu_{0}, u) = \lambda \cdot f(x, u) - L(x, u),
\end{equation*}
and the optimal contact Hamiltonian~\eqref{eq:h} becomes
\begin{equation}
  \label{eq:h-normal}
  h(\hat{x}, [\hat{\nu}]) = \lambda \cdot f(x, u^{\star}(x,\lambda)) - L(x, u^{\star}(x,\lambda)).
\end{equation}
Therefore, {\em the conventional practice of getting rid of the redundancy in the costate vector $\hat{\nu} \in \R^{n+1}$ for normal extremals by setting $\nu_{0} = -1$ is equivalent to regarding the costate as an element $[\hat{\nu}]$ in the projective space $\mathbb{P}(\R^{n+1})$ and expressing it in terms of appropriate homogeneous coordinates for $\mathbb{P}(\R^{n+1})$.}

\begin{figure}[htbp]
  \centering
  \includegraphics[width=.35\linewidth]{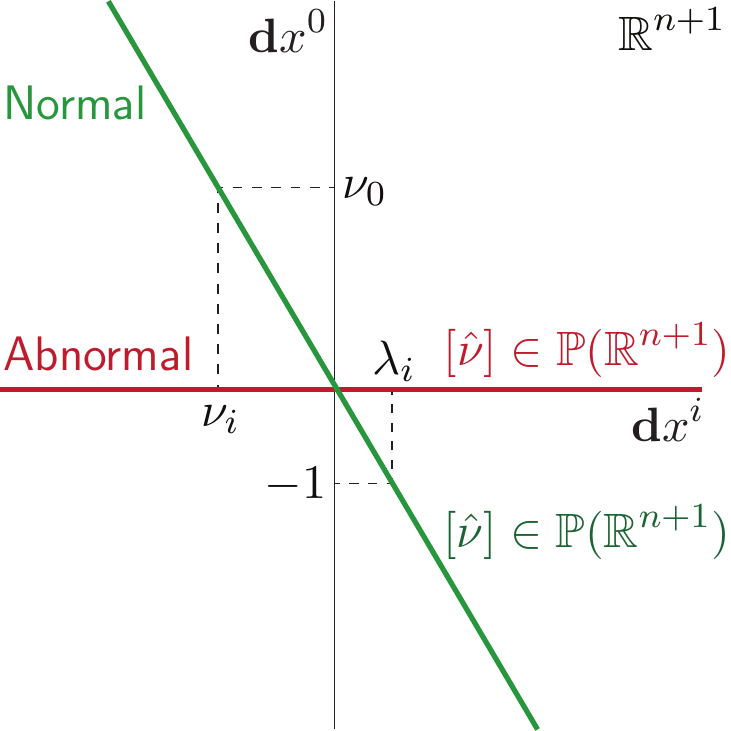}
  \caption{The costates for normal and abnormal minimizers.}
  \label{fig:Normal-Abnormal}
\end{figure}

Those lines or costates $[\hat{\nu}]$ with $\nu_{0} = 0$ are at the {\em coordinate singularity} of the above homogeneous coordinates $\lambda$ and so one needs to employ a different coordinate chart for such costate $[\hat{\nu}]$:
\begin{equation*}
  [\hat{\nu}] = \brackets{ 0: \alpha_{1} : \dots : \alpha_{n}},
\end{equation*}
where $\alpha_{i} = \nu_{i}/\nu_{a}$ for some fixed $a \in \{1, \dots, n\}$ such that $\nu_{a} \neq 0$.
So we may now define the contact control Hamiltonian~\eqref{eq:h_c} as
\begin{equation*}
  h_{\rm c}(\hat{x}, [\hat{\nu}], u) \defeq H_{\rm c}(\hat{x}, \hat{\nu}/\nu_{a}, u) = \alpha \cdot f(x, u).
\end{equation*}
Therefore, {\em an abnormal extremal is identified as a costate $[\hat{\nu}]$ at the coordinate singularity of the standard homogeneous coordinates~\eqref{eq:homogeneous_coordinates} (see Fig~\ref{fig:Normal-Abnormal})}; so the use of the projective space $\mathbb{P}(\R^{n+1})$ gives rise to a differential-geometric classification of normal and abnormal extremals.

\subsection{Hyperplane Field and Manifold of Contact Elements}
\label{ssec:HyperplaneField}
Recall from Section~\ref{ssec:Costate_and_ProjectiveSpace} that we introduced the following identification (see Fig.~\ref{fig:HyperplaneAndRay}):
\begin{equation*}
  \text{separating hyperplane $\mathcal{H}_{\hat{x}_{1}^{\star}}$ in $\R^{n+1}$}
  \leftrightarrow
  \text{costate $[\hat{\nu}^{\star}(t_{1}^{\star})]$ in $\mathbb{P}(\R^{n+1})$ s.t. $\mathcal{H}_{\hat{x}_{1}^{\star}} = \ker\hat{\nu}^{\star}(t_{1}^{\star})$}.
\end{equation*}
In fact, {\em the idea of identifying a hyperplane in $\R^{n+1}$ with an element in the projective space $\mathbb{P}(\R^{n+1})$ is essential in contact geometry.}
One may define the set $\mathbb{H}$ of hyperplanes in $\R^{n+1}$, where we identify those hyperplanes that are translations of one another and so a single representative element in $\mathbb{H}$ would be a hyperplane $\mathcal{H}$ passing through the origin.
Then one easily sees that $\mathbb{H}$ is identified with $\mathbb{P}(\R^{n+1})$; therefore,
\begin{equation*}
  \mathbb{H} \defeq \braces{ \text{hyperplanes in $\R^{n+1}$} }
  \cong
  \braces{ \text{costates} } = \mathbb{P}(\R^{n+1}).
\end{equation*}

\begin{figure}[htbp]
  \centering
  \includegraphics[width=.5\linewidth]{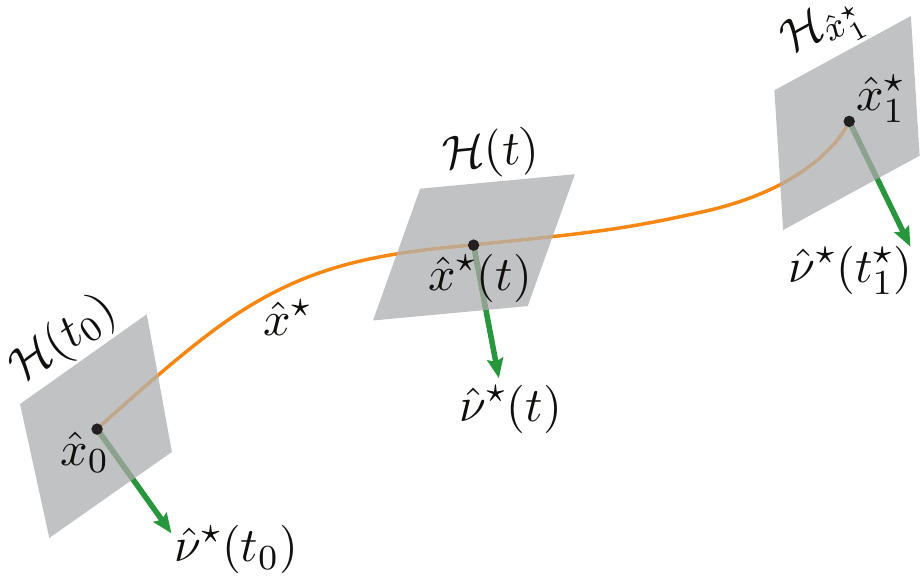}
  \caption{Hyperplane field propagated along optimal solution.}
  \label{fig:HyperplaneField}
\end{figure}

Now, recall that we employed the adjoint equation~\eqref{eq:adjoint} to propagate the costate vector $\hat{\nu}^{\star}(t)$ along the optimal solution $\hat{x}^{\star}(t)$.
We may now see the pair $(\hat{x}^{\star}(t), [\hat{\nu}^{\star}(t)])$ as a curve in $\R^{n+1} \times \mathbb{P}(\R^{n+1})$; alternatively, we may define the {\em hyperplane field} $\mathcal{H}(t) \defeq \ker\hat{\nu}^{\star}(t)$ along $\hat{x}^{\star}(t)$ and see the pair $(\hat{x}^{\star}(t), \mathcal{H}(t))$ as a curve in $\R^{n+1} \times \mathbb{H}$ (see Fig.~\ref{fig:HyperplaneField}).
They are two different pictures of the same thing.
In fact, the identification $\R^{n+1} \times \mathbb{H} \cong \R^{n+1} \times \mathbb{P}(\R^{n+1})$ is standard in contact geometry, and they are a simple example of {\em manifold of contact elements} and is a basic example of {\em contact manifold} as well.

\subsection{Digression on Contact Geometry}
\label{sec:DigressionOnContactGeom}
Saving a more general treatment for later (see Section~\ref{ssec:ContactGeometry}), this subsection briefly explains what makes the spaces $\R^{n+1} \times \mathbb{H}$ and $\R^{n+1} \times \mathbb{P}(\R^{n+1})$ introduced above contact manifolds.

First we define some shorthand notation:
\begin{equation*}
  \mathcal{C} \defeq \R^{n+1} \times \mathbb{H},
  \qquad
  \mathbb{P}(T^{*}\R^{n+1}) \defeq \R^{n+1} \times \mathbb{P}(\R^{n+1}).
\end{equation*}
Since the space $\mathbb{H}$ is the collection of hyperplanes (passing through the origin) in $\R^{n+1}$, each element in $\mathcal{C}$ may be regarded as the pair $\mathcal{H}_{\hat{x}} \defeq (\hat{x}, \mathcal{H})$ of a base point $\hat{x} \in \R^{n+1}$ and a hyperplane $\mathcal{H} \subset \R^{n+1}$ attached to $\hat{x}$; thus the space $\mathcal{C}$ is the collection of all hyperplane fields on $\R^{n+1}$.
Likewise, each element $[\hat{\nu}]_{\hat{x}} \defeq (\hat{x}, [\hat{\nu}])$ in $\mathbb{P}(T^{*}\R^{n+1})$ may be regarded as an assignment of a straight line in $\R^{n+1}$ to a base point $\hat{x} \in \R^{n+1}$.
We may identify $\mathcal{C}$ with $\mathbb{P}(T^{*}\R^{n+1})$ by the relation $\mathcal{H} = \ker\hat{\nu}$ (recall Fig.~\ref{fig:HyperplaneAndRay}).
Using the homogeneous coordinates $\lambda \defeq (\lambda_{1}, \dots, \lambda_{n})$ for $\mathbb{P}(\R^{n+1})$ from \eqref{eq:homogeneous_coordinates}, we assign coordinates $(x^{0}, x, \lambda)$ to an element $(\hat{x}, \mathcal{H}) \in \mathcal{C}$ or $(\hat{x}, [\hat{\nu}]) \in \mathbb{P}(T^{*}\R^{n+1})$; thus $\mathcal{C} \cong \mathbb{P}(T^{*}\R^{n+1})$ is a $(2n+1)$-dimensional manifold.

What makes $\mathcal{C} \cong \mathbb{P}(T^{*}\R^{n+1})$ a {\em contact manifold} is the one-form $\theta$ on it called a {\em contact form} that is written as
\begin{equation*}
  \theta_{[\hat{\nu}]_{\hat{x}}} = -\d{x^{0}} + \lambda_{i}\,\d{x^{i}},
\end{equation*}
where $\d$ is the exterior differential.
The defining characteristics of the contact form $\theta$ is that it defines a hyperplane $\ker\theta_{[\hat{\nu}]_{\hat{x}}}$ at each point $[\hat{\nu}]_{\hat{x}}$ on $\mathcal{C}$ so that the two-form
\begin{equation*}
  \omega_{[\hat{\nu}]_{\hat{x}}} \defeq -\d\theta_{[\hat{\nu}]_{\hat{x}}} = \d{x^{i}} \wedge \d\lambda_{i}
\end{equation*}
is non-degenerate on the hyperplane $\ker\theta_{[\hat{\nu}]_{\hat{x}}}$; such a hyperplane assignment is called a {\em contact structure} which, by definition, makes $\mathcal{C}$ a {\em contact manifold}; see, e.g., \citet[Appendix~4]{Ar1991}, \citet[Chapter~10]{KuLyRu2007}, \citet[Chapters~10 \& 11]{Ca2008}, and \citet[Chapters~1 \& 2]{Ge2008} for details.

\subsection{Contact Hamiltonian System and Necessary Condition for Optimality}
\label{ssec:ContactHamSysAndNecessaryCond}
Given a function (contact Hamiltonian) $h\colon \mathcal{C} \cong \mathbb{P}(T^{*}\R^{n+1}) \to \R$, one may define the corresponding {\em contact Hamiltonian vector field} $X_{h}$ on $\mathcal{C} \cong \mathbb{P}(T^{*}\R^{n+1})$ as follows:
\begin{equation}
\label{eq:ContactFlow}
  \theta(X_{h}) = h,
  \qquad
  \ins{X_{h}} \omega = \d{}h - (\d{}h \cdot R_{\theta})\,\theta,
\end{equation}
where $R_{\theta}$ is the {\em Reeb vector field} associated with the contact form $\theta$, i.e., the vector field $R_{\theta}$ on $\mathcal{C}$ that is uniquely defined by
\begin{equation*}
  \theta(R_{\theta}) = 1,
  \qquad
  \ins{R_{\theta}} \omega = 0,
\end{equation*}
which gives
\begin{equation*}
  R_{\theta}(\hat{x}, \lambda) = -\pd{}{x^{0}} = -\hat{\bf e}_{0}.
\end{equation*}
Hence the Reeb vector field $R_{\theta}$ defines one of the ``forbidden'' directions in $\R_{\le0} \times T_{x^{\star}_1}S_{1}$ (see Fig.~\ref{fig:SeparatingHyperplane}).
For a normal extremal, we may use the coordinates $(x^{0}, x, \lambda)$ to write \eqref{eq:ContactFlow} as follows:
\begin{equation*}
  \dot{x}^{0} = h - \lambda_{i} \dot{x}^{i},
  \qquad
  \dot{x}^{i} = \pd{h}{\lambda_{i}},
  \qquad
  \dot{\lambda}_{i} = \lambda_{i} \pd{h}{x^{0}} - \pd{h}{x^{i}}.
\end{equation*}
In particular, with the optimal contact Hamiltonian~\eqref{eq:h-normal}, we have
\begin{equation*}
  \dot{x}^{\star0} = L(x^{\star}, u^{\star}),
  \qquad
  \dot{x}^{\star i} = f^{i}(x^{\star}, u^{\star}),
  \qquad
  \dot{\lambda}^{\star}_{i} = -\lambda^{\star}_{j} \pd{f^{j}}{x^{i}}(x^{\star}, u^{\star}) + \pd{L}{x^{i}}(x^{\star}, u^{\star}),
\end{equation*}
which is the extended system~\eqref{eq:extended_system} along with the adjoint equation~\eqref{eq:adjoint} for $u = u^{\star}$.
A similar result follows for an abnormal extremal as well.
To summarize, we have the following:
\begin{proposition}
  \label{prop:ContactHamSysAndNecessaryCond}
  Let $u^{\star}\colon [t_{0}, t_{1}^{\star}] \to \mathcal{U}$ be an optimal control.
  Then the corresponding optimal trajectory and costate $(\hat{x}^{\star}(t), [\hat{\nu}^{\star}](t))$ for $t \in [t_{0}, t_{1}^{\star}]$ satisfy the contact Hamiltonian system~\eqref{eq:ContactFlow} corresponding to the optimal contact Hamiltonian $h(\cdot, u^{\star}(\cdot))\colon \mathbb{P}(T^{*}\R^{n+1}) \to \R$.
\end{proposition}
\begin{remark}
  \label{rem:symplectic_vs_contact}
  One may formulate the necessary condition on the cotangent bundle $T^{*}\R^{n+1}$ as a Hamiltonian system in the symplectic sense.
  However, this geometric setting is less suited for interpreting the duality between the separating hyperplanes and the costate (discussed in Section~\ref{ssec:HyperplaneField}).
  The symplectic formulation involves the coordinate $\nu_{0}$ for the costate vector $\hat{\nu} \in T^{*}\R^{n+1}$, but this is redundant because the corresponding Hamilton equations give $\dot{\nu}_{0} = 0$; see \eqref{eq:adjoint}.
  The contact-geometric view gets rid of this redundancy at the outset by projectivization; as a result, it naturally gives rise to the identification of hyperplanes and projective cotangent spaces.
  This is exactly the duality between the separating hyperplanes and the costate exploited in the maximum principle, as explained in Sections~\ref{ssec:Costate_and_ProjectiveSpace} and \ref{ssec:HyperplaneField}.
\end{remark}

\section{Geometry of Optimal Control on Manifolds}
We now replace the state space $\R^{n}$ by an $n$-dimensional manifold $M$ to consider optimal control of systems on the manifold $M$.
We do not delve into the details on how to extend the proof of the maximum principle to manifolds; see, e.g., \citet[Chapter~11]{Ju1997}, \citet[Chapter~12]{AgSa2004}, \citet{BaMu2009}, and \citet{Ch2011}.
Skipping technical details, we briefly sketch how the geometric ideas in $\R^{n}$ from the previous section can be generalized to manifolds.
Our main references in contact geometry are \citet[Appendix~4]{Ar1991}, \citet[Chapter~10]{KuLyRu2007}, \citet[Chapters~10 \& 11]{Ca2008}, and \citet[Chapters~1 \& 2]{Ge2008}.

\subsection{Optimal Control on Manifolds}
Let $\tau\colon TM \to M$ be the tangent bundle of $M$, $\mathcal{U}$ a compact subset of $\R^{m}$ as before, and $pr\colon M \times \mathcal{U} \to M$ the projection to the first slot.
We now have a map $f\colon M \times \mathcal{U} \to TM$ such that $\tau \circ f = pr$ and a cost function $L\colon M \times \mathcal{U} \to \R$, and then we can formulate the optimal control problem on $M$ just as in \eqref{eq:OptCtrlProblem}.

Then we define the extended configuration space $\hat{M} \defeq \R \times M = \{\hat{x} \defeq (x^{0}, x)\}$, and also $\hat{f}\colon \hat{M} \times \mathcal{U} \to \R \times TM$ by $\hat{f}( \hat{x}, u) \defeq \parentheses{ L(x,u), f(x,u) }$; hence we may define the extended system~\eqref{eq:extended_system}.

The cone $C_{\hat{x}_{1}^{\star}}$ is now more naturally a subset of the tangent space $T_{\hat{x}_{1}^{\star}}\hat{M}$ and the costate vector $\hat{\nu}^{\star}(t_{1}^{\star})$ is in the cotangent space $T_{\hat{x}_{1}^{\star}}^{*}\hat{M}$.
Likewise, the hyperplane $\mathcal{H}(t)$ and costate vector $\hat{\nu}^{\star}(t)$ at $\hat{x}^{\star}(t)$ are in the tangent and cotangent spaces, respectively, i.e., $\mathcal{H}(t) \subset T_{\hat{x}^{\star}(t)}\hat{M}$ and $\hat{\nu}^{\star}(t) \in T_{\hat{x}^{\star}(t)}^{*}\hat{M}$.
As briefly mentioned in Remark~\ref{rem:cotangent_lift}, one needs to propagate back the costate vector $\hat{\nu}^{\star}(t_{1}^{\star})$ along $\hat{x}^{\star}(t)$ by the cotangent lift of the flow $\phi_{t_{1}^{\star}-t}\colon \hat{M} \to \hat{M}$ defined by the optimal solution $\dot{\hat{x}}^{\star} = \hat{f}(\hat{x}^{\star}, u^{\star})$, i.e., $\hat{\nu}^{\star}(t) \defeq T^{*}\phi_{t_{1}^{\star}-t}(\hat{\nu}^{\star}(t_{1}^{\star}))$, and hence
\begin{equation*}
  [\hat{\nu}^{\star}(t)] =  [T^{*}\phi_{t_{1}^{\star}-t}(\hat{\nu}^{\star}(t_{1}^{\star}))] \in \mathbb{P}(T_{\hat{x}^{\star}(t)}\hat{M}),
\end{equation*}
where $\mathbb{P}(T_{\hat{x}^{\star}(t)}\hat{M})$ is the projectivization of the cotangent space $T_{\hat{x}^{\star}(t)}^{*}\hat{M}$.

\subsection{Contact Geometry and Necessary Condition}
\label{ssec:ContactGeometry}
The space $\mathcal{C} \defeq \R^{n+1} \times \mathbb{H}$ introduced in Section~\ref{sec:DigressionOnContactGeom} is our prototype of what is called a {\em manifold of contact elements}, which we define now:
Let $\hat{M}$ be a manifold and $T\hat{M}$ its tangent bundle.
A {\em contact element} on $\hat{M}$ is a point $\hat{x} \in \hat{M}$ along with a hyperplane (passing through the origin) $\mathcal{H}_{\hat{x}} \subset T_{\hat{x}}\hat{M}$.
The collection $\mathcal{C}$ of contact elements $(\hat{x}, \mathcal{H}_{\hat{x}})$ is called a {\em manifold of contact elements} of $\hat{M}$.
Note that an element in the manifold $\mathcal{C}$ is a point $\hat{x} \in \hat{M}$ along with a hyperplane $\mathcal{H}_{\hat{x}} \subset T_{\hat{x}}\hat{M}$ attached to the point $\hat{x}$.

Likewise, the space $\mathbb{P}(T^{*}\R^{n+1}) \defeq \R^{n+1} \times \mathbb{P}(\R^{n+1})$ from Section~\ref{sec:DigressionOnContactGeom} is an example of the {\em projectivized cotangent bundle} $\mathbb{P}(T^{*}\hat{M})$ defined by
\begin{equation*}
   \mathbb{P}(T^{*}\hat{M}) \defeq \bigcup_{\hat{x} \in \hat{M}} \mathbb{P}(T^{*}_{\hat{x}}\hat{M}).
\end{equation*}
Then the hyperplane $\mathcal{H}_{\hat{x}} \subset T_{\hat{x}}\hat{M}$ is identified with $[\hat{\nu}]_{\hat{x}} \in \mathbb{P}(T_{\hat{x}}^{*}\hat{M})$, and thus $\mathcal{C}$ is identified with the projectivized cotangent bundle $\mathbb{P}(T^{*}\hat{M})$.
Note that then $[\hat{\nu}^{\star}(t)]$ and $\mathcal{H}(t) \defeq \ker\hat{\nu}^{\star}(t)$ are curves in $\mathbb{P}(T^{*}\hat{M})$ and $\mathcal{C}$, respectively.
An element $[\hat{\nu}]_{\hat{x}}$ in each fiber $\mathbb{P}(T^{*}_{\hat{x}}\hat{M})$ is parametrized by the homogeneous coordinates $\{\lambda_{i} \defeq -\nu_{i}/\nu_{0}\}_{i=1}^{n}$ defined in \eqref{eq:homogeneous_coordinates} for normal extremals.
Therefore, $(x^{0}, \dots, x^{n}, \lambda_{1}, \dots, \lambda_{n})$ gives local coordinates for $\mathbb{P}(T^{*}\hat{M})$ just as with the case with $\R^{n+1} \times \mathbb{P}(\R^{n+1})$; in fact, for $\hat{M} = \R^{n+1}$, $\mathbb{P}(T^{*}\hat{M}) \cong \R^{n+1} \times \mathbb{P}(\R^{n+1})$.

Let $\pi\colon T^{*}\hat{M} \to \hat{M}$ be the cotangent bundle projection and $[\pi]\colon \mathbb{P}(T^{*}\hat{M}) \to \hat{M}$ be its projectivization, i.e., $[\pi]([\hat{\nu}]_{\hat{x}}) = \hat{x}$ for any $[\hat{\nu}]_{\hat{x}} \defeq (\hat{x}, [\hat{\nu}]) \in \mathbb{P}(T^{*}\hat{M})$.
Now, we define a one-form $\theta$ on $\mathbb{P}(T^{*}\hat{M})$ as follows:
\begin{equation*}
  \theta_{[\hat{\nu}]_{\hat{x}}} \cdot w_{[\hat{\nu}]_{\hat{x}}} = \hat{\nu}_{\hat{x}} \cdot T_{[\hat{\nu}]_{\hat{x}}}[\pi](w_{[\hat{\nu}]_{\hat{x}}})
\end{equation*}
for any $w_{[\hat{\nu}]_{\hat{x}}} \in T_{[\hat{\nu}]_{\hat{x}}}\mathbb{P}(T^{*}\hat{M})$.
We also define a two-form $\omega$ on $\mathbb{P}(T^{*}\hat{M})$ by $\omega \defeq -\d\theta$.
Locally, we have
\begin{equation*}
  \theta_{[\hat{\nu}]_{\hat{x}}} = -\d{x^{0}} + \lambda_{1}\d{x^{1}} + \dots + \lambda_{n}\d{x^{n}}
\end{equation*}
and
\begin{equation*}
  \omega_{[\hat{\nu}]_{\hat{x}}} = \d{x^{1}} \wedge \d\lambda_{1} + \dots + \d{x^{n}} \wedge \d\lambda_{n}.
\end{equation*}
The one-form $\theta$ then defines the hyperplane $\ker\theta_{[\hat{\nu}]_{\hat{x}}}$ at every point $[\hat{\nu}]_{\hat{x}}$ of $\mathbb{P}(T^{*}\hat{M})$ so that $\omega$ is non-degenerate on the hyperplane, i.e., the hyperplane field $\ker\theta$ defines a {\em contact structure} on $\mathbb{P}(T^{*}\hat{M})$, which makes itself a {\em contact manifold}.

We may then define contact Hamiltonian systems just as in Section~\ref{ssec:ContactHamSysAndNecessaryCond} and can generalize Proposition~\ref{prop:ContactHamSysAndNecessaryCond} to control systems on the manifold $M$, where the optimal contact Hamiltonian $h$ is defined on $\mathbb{P}(T^{*}\hat{M})$.

\section{Application:  Terminal Cost and Transversality Condition}
\subsection{Optimal Control to a General Target with Terminal Cost}
Let $K\colon \R^{n} \to \R$ be a smooth function defined on the configuration space $M = \R^{n}$ and consider the following variant of the optimal control problem~\eqref{eq:OptCtrlProblem} with the terminal cost $K(x(t_{1}))$:
\begin{equation*}
  \begin{array}{c}
    \DS \min_{u(\cdot) \in \mathcal{U}} \brackets{ K(x(t_{1})) + \int_{t_{0}}^{t_{1}} L(x(t), u(t))\,dt }
    \medskip\\
    \text{subject to}
    \quad
    \dot{x} = f(x, u),
    \quad
    x(t_{0}) = x_{0},
    \quad
    x(t_{1}) \in S_{1},
  \end{array}
\end{equation*}
where the conditions for the end points are the same as before.
One may now define
\begin{equation*}
  x^{0}(t) \defeq K(x(t)) + \int_{t_{0}}^{t} L(x(s), u(s))\,ds
\end{equation*}
with the initial condition $x^{0}(t_{0}) = K(x(t_{0}))$.
Then we have
\begin{equation*}
  \od{}{t}[ x^{0}(t) - K(x(t)) ] = L(x(t), u(t)),
\end{equation*}
or defining $\hat{y} = (y^{0}, y) \defeq (x^{0} - K(x), x)$, we have
\begin{equation}
  \label{eq:extended_system-y}
  \dot{y}^{0} = L(y, u),
  \qquad
  \dot{y} = f(y, u)
\end{equation}
with the initial condition $y^{0}(t_{0}) = 0$; thus we have the same extended system as in \eqref{eq:extended_system}.
Therefore, we may apply Proposition~\ref{prop:ContactHamSysAndNecessaryCond} to the extended system \eqref{eq:extended_system-y}, and so the optimal flow is given by the contact Hamiltonian flow $X_{h}$ with the contact Hamiltonian $h\colon \mathbb{P}(T^{*}\R^{n+1}) \to \R$ defined by
\begin{equation*}
  h(\hat{y}, \mu) \defeq \mu_{i} f^{i}(y, u^{\star}) - L(y, u^{\star}),
\end{equation*}
where $[(-1, \mu)] \in \mathbb{P}(T^{*}\R^{n+1})$ is the costate corresponding to $y$.

\subsection{Transversality Condition via Contact Transformation}
One may be then tempted to write $\mu^{\star}(t_{1}^{\star}) \in ( T_{x^{\star}_{1}}S_{1} )^{\perp}$ as the transversality condition as in \eqref{eq:transversality_condition}, but this is incorrect because $[(-1, \mu)] \in \mathbb{P}(T^{*}\R^{n+1})$ is the costate corresponding to the state variables $(y^{0}, y) \in \R^{n+1}$ and
\begin{equation*}
  y^{0}(t_{1}) = \int_{t_{0}}^{t_{1}} L(x(t), u(t))\,dt 
\end{equation*}
is {\em not} the quantity to be minimized.
Instead, it is $x^{0}(t_{1})$ that is to be minimized, and therefore we may write the {\em correct} transversality condition
\begin{equation}
  \label{eq:transversality_condition-lambda}
  \lambda^{\star}(t_{1}^{\star}) \in (T_{x^{\star}_{1}}S_{1})^{\perp}  
\end{equation}
for the costate $[(-1,\lambda)] \in \mathbb{P}(T^{*}\R^{n+1})$ corresponding to the original variable $x$.

The question is then: How does one rewrite \eqref{eq:transversality_condition-lambda} in terms of $\mu$?
Contact geometry provides a simple and elegant answer to this question and leads us to a simple derivation of the correct transversality condition for $\mu$:
The discussion in the previous subsection motivates us to define the diffeomorphism $\Phi_{K}\colon \R^{n+1} \to \R^{n+1}$ defined by
\begin{equation*}
  \Phi_{K}(x^{0}, x) \defeq \parentheses{ x^{0} - K(x), x } = (y^{0}, y).
\end{equation*}
Clearly its inverse is given by $\Phi_{K}^{-1} = \Phi_{-K}$.
Now let $\mathcal{H}_{\hat{x}_{1}^{\star}}$ be the separating hyperplane at $\hat{x}_{1}^{\star}$.
Then $T_{\hat{x}_{1}^{\star}}\Phi_{K}(\mathcal{H}_{\hat{x}_{1}^{\star}})$ gives the separating hyperplane at $\hat{y}_{1}^{\star}$, and we have
\begin{equation*}
  0 = \ip{ \hat{\nu}_{\hat{x}_{1}^{\star}} }{ \mathcal{H}_{\hat{x}_{1}^{\star}} } = \ip{ T^{*}_{\hat{x}_{1}^{\star}}\Phi_{K}^{-1}(\hat{\nu}_{\hat{x}_{1}^{\star}}) }{ T_{\hat{x}_{1}^{\star}}\Phi_{K}(\mathcal{H}_{\hat{x}_{1}^{\star}}) }.
\end{equation*}
So $T^{*}_{\hat{x}_{1}^{\star}}\Phi_{K}^{-1}(\hat{\nu}_{\hat{x}_{1}^{\star}})$ is a costate vector in $T_{\hat{y}_{1}^{\star}}^{*}\R^{n+1}$ and hence
\begin{equation*}
  [(-1,\mu^{\star}(t_{1}^{\star}))] = [T^{*}_{\hat{x}_{1}^{\star}}\Phi_{K}^{-1}(\hat{\nu}_{\hat{x}_{1}^{\star}})],
\end{equation*}
where $[\hat{\nu}_{\hat{x}_{1}^{\star}}] = [(-1,\lambda^{\star}(t_{1}^{\star}))]$.
As a result, the costates $[(-1,\lambda)]$ and $[(-1,\mu)]$ are related by the projectivization
\begin{equation*}
  \Psi_{K} \defeq [T^{*}\Phi_{K}]\colon \mathbb{P}(T^{*}\R^{n+1}) \to \mathbb{P}(T^{*}\R^{n+1})
\end{equation*}
of $T^{*}\Phi_{K}\colon T^{*}\R^{n+1} \to T^{*}\R^{n+1}$; specifically,
\begin{equation*}
  (\hat{x}, \lambda) = \Psi_{K}(\hat{y}, \mu),
\end{equation*}
i.e., the diagram below commutes.
\begin{equation*}
  \begin{tikzcd}[column sep=5.5ex, row sep=6.5ex]
    \mathbb{P}(T^{*}\R^{n+1}) \arrow{d}[swap]{[\pi]} & \mathbb{P}(T^{*}\R^{n+1}) \arrow{d}{[\pi]} \arrow{l}[swap]{\Psi_{K}}
    \\
    \R^{n+1} \arrow{r}[swap]{\Phi_{K}} & \R^{n+1}
  \end{tikzcd}
  \quad
  \begin{tikzcd}[column sep=6ex, row sep=6.5ex]
    (\hat{x},\lambda) \arrow[mapsto]{d} & (\hat{y},\mu) \arrow[mapsto]{d} \arrow[mapsto]{l}
    \\
    \hat{x} \arrow[mapsto]{r} & \hat{y}
  \end{tikzcd}
\end{equation*}
Therefore, we have
\begin{equation*}
  (x^{0}, x, \lambda) = \Psi_{K}(y^{0}, y, \mu) = \parentheses{ y^{0} + K(y),\; y,\; \mu + \d{K}(y) },
\end{equation*}
and then applying \eqref{eq:transversality_condition-lambda} to the above equation gives the {\em correct} transversality condition for $\mu$:
\begin{equation*}
  \mu^{\star}(t_{1}^{\star}) + \d{K}(x^{\star}_{1}) \in ( T_{x^{\star}_{1}}S_{1} )^{\perp}.
\end{equation*}

\begin{remark}
  One may, in principle, work with costate vectors in $T^{*}\R^{n+1}$ without projectivization, but as mentioned in Remark~\ref{rem:symplectic_vs_contact}, the result comes with a redundancy in the costate vector, i.e., $\mu_{0} = \lambda_{0}$ but $\lambda_{0}$ is left arbitrary; whereas the projectivization gets rid of the redundancy at the outset.
\end{remark}

\bibliography{ContactGeomOfPMP}
\bibliographystyle{plainnat}

\end{document}